\documentstyle[11pt]{article} 
\begin{document}
\def\r{\rightarrow}
\newcommand{\fdem}{\hfill$\Box$}

\newcommand{\E}{{I\!\!E}}     
\renewcommand{\P}{{I\!\!P}}     
\renewcommand{\L}{{I\!\!L}}     
\newcommand{\F}{{I\!\!F}}     
\newcommand{\T}{\mbox{\rule{0.1ex}{1.5ex}\hspace{0.2ex}\rule{0.1ex}{1.5ex}
\hspace{-1.7ex}\rule[1.5ex]{1.6ex}{0.1ex} }}
     
\newcommand{\N}{{I\!\!N}}     
\newcommand{\Z}{{Z\!\!\!Z}}
\newcommand{\R}{{I\!\!R}}     
\newcommand{\D}{{I\!\!D}}     
\newcommand{\C}{\mbox{\rm C\hspace{-1.1ex}\rule[0.3ex]{0.1ex}{1.2ex}
\hspace{1.1ex}}}
\newcommand{\Q}{\mbox{\rm Q\hspace{-1.1ex}\rule[0.3ex]{0.1ex}{1.2ex}
\hspace{1.1ex}}}

\renewcommand{\Im}{\mbox{\rm Im }}
\renewcommand{\Re}{\mbox{\rm Re }}
\newcommand{\supp}{\mathop{\rm supp}}
\newcommand{\diag}{\mathop{\rm diag}}
\renewcommand{\deg}{\mathop{\rm deg}}
\renewcommand{\dim}{\mathop{\rm dim}}
\renewcommand{\ker}{\mathop{\rm Ker}}
\newcommand{\id}{\mbox{\rm1\hspace{-.2ex}\rule{.1ex}{1.44ex}}
   \hspace{-.82ex}\rule[-.01ex]{1.07ex}{.1ex}\hspace{.2ex}}
\newcommand{\bra}{\langle\,}
\newcommand{\ket}{\,\rangle}
\newcommand{\obl}{/\!/}
\newcommand{\mapdown}[1]{\vbox{\vskip 4.25pt\hbox{\bigg\downarrow
  \rlap{$\vcenter{\hbox{$#1$}}$}}\vskip 1pt}}
\newcommand{\tr}{\mbox{\rm tr}}
\newcommand{\Tr}{\mbox{\rm Tr}}
\newcommand{\r}{\mathop{\rightarrow}}
\newcommand{\cA}{\mbox{$\cal A$}}
\newcommand{\cB}{\mbox{$\cal B$}}
\newcommand{\cC}{\mbox{$\cal C$}}
\newcommand{\cD}{\mbox{$\cal D$}}
\newcommand{\cE}{\mbox{$\cal E$}}
\newcommand{\cF}{\mbox{$\cal F$}}
\newcommand{\cG}{\mbox{$\cal G$}}
\newcommand{\cH}{\mbox{$\cal H$}}
\newcommand{\cJ}{\mbox{$\cal J$}}
\newcommand{\cK}{\mbox{$\cal K$}}
\newcommand{\cL}{\mbox{$\cal L$}}
\newcommand{\cM}{\mbox{$\cal M$}}
\newcommand{\cN}{\mbox{$\cal N$}}
\newcommand{\cP}{\mbox{$\cal P$}}
\newcommand{\cR}{\mbox{$\cal R$}}
\newcommand{\cS}{\mbox{$\cal S$}}
\newcommand{\cT}{\mbox{$\cal T$}}
\newcommand{\cW}{\mbox{$\cal W$}}
\newcommand{\cY}{\mbox{$\cal Y$}}
\newcommand{\cX}{\mbox{$\cal X$}}
\newcommand{\cU}{\mbox{$\cal U$}}
\newcommand{\cV}{\mbox{$\cal V$}}

\begin{center}
{\large\bf VITESSE DE CONVERGENCE DANS LE TH\'EOR\`EME 

LIMITE CENTRAL POUR CHA\^INES DE MARKOV 

DE PROBABILIT\'E DE TRANSITION QUASI-COMPACTE.}

\end{center}
\vskip 7mm
\begin{center}
HERV\'E Lo\"{\i}c \\
I.R.M.A.R., UMR-CNRS 6625, \\ 
Institut National des Sciences Appliqu\'ees de Rennes, \\ 
20, Avenue des Buttes de Cou\"esmes CS 14 315, 35043 Rennes Cedex. \\ 
Loic.Herve@insa-rennes.fr
\end{center}




\vskip 7mm

\noindent{\bf R\'esum\'e.} {\small \it Soit $Q$ une probabilit\'e de transition sur un espace mesurable $E$, 
admettant une probabilit\'e invariante, 
soit $(X_n)_n$ une cha\^{\i}ne de Markov associ\'ee  \`a $Q$, et enfin soit $\xi$ une fonction r\'eelle mesurable sur $E$, et 
$S_n = \sum_{k=1}^{n} \xi(X_k)$. 
Sous des hypoth\`eses fonctionnelles sur l'action de $Q$ et de ses noyaux de Fourier $Q(t)$, nous \'etudions la vitesse de convergence
dans le th\'eor\`eme limite central pour la suite $(\frac{S_n}{\sqrt n})_n$. Selon les hypoth\`eses nous obtenons une vitesse 
en $n^{-\frac{\tau}{2}}$ pour tout $\tau<1$, ou bien en $n^{-\frac{1}{2}}$. 
Nous appliquons la m\'ethode spectrale de Nagaev en l'am\'eliorant, d'une part gr\^ace \`a un 
th\'eor\`eme de perturbations de Keller et Liverani, d'autre part gr\^ace \`a une majoration de 
$\E[e^{it\frac{S_n}{\sqrt n}}] - e^{\frac{-t^2}{2}}$ obtenue par une m\'ethode de r\'eduction en diff\'erence de martingale. 
Lorsque $E$ est non compact ou $\xi$ est non born\'ee, 
les conditions requises ici sur $Q(t)$ (en substance, des conditions de moment sur $\xi$) sont plus faibles que celles 
habituellement impos\'ees lorsqu'on utilise le th\'eor\`eme de perturbation standard. Par exemple, dans le cadre 
des cha\^{\i}nes $V$-g\'eom\'etriquement ergodiques ou des mod\`eles it\'eratifs Lipschitziens, on obtient dans le t.l.c 
une vitesse en $n^{-\frac{1}{2}}$ sous une hypoth\`ese de moment d'ordre 3 sur $\xi$. }

\vskip 5mm 

\noindent{\bf Abstract.} {\small \it 
Let $Q$ be a transition probability on a measurable space $E$ which admits an invariant probability measure, 
let $(X_n)_n$ be a Markov chain associated to $Q$, and let 
$\xi$ be a real-valued measurable function on $E$, and $S_n = \sum_{k=1}^{n} \xi(X_k)$. 
Under functional hypotheses on the action of $Q$ and its Fourier kernels $Q(t)$, 
we investigate the rate of convergence in the central limit theorem for the sequence 
$(\frac{S_n}{\sqrt n})_n$. According to the hypotheses, we prove that the rate is, either $O(n^{-\frac{\tau}{2}})$ for all $\tau<1$, or 
$O(n^{-\frac{1}{2}})$.
We apply the spectral method of Nagaev which is improved by using a perturbation theorem of Keller and Liverani, and 
a majoration of $\E[e^{it\frac{S_n}{\sqrt n}}] - e^{\frac{-t^2}{2}}$ obtained by a method of martingale difference reduction.  
When $E$ is not compact or $\xi$ is not bounded, the conditions required here on $Q(t)$ (in substance, some moment conditions on $\xi$) 
are weaker than the ones usually imposed when the standard perturbation theorem is used in the spectral method. 
For example, in the case of $V$-geometric ergodic chains
or Lipschitz iterative models, the rate of convergence in the c.l.t is $O(n^{-\frac{1}{2}})$ under a third moment condition on $\xi$. }


\newpage 

\pagestyle{myheadings}
\markright{\small $\ \ \ \ \ \ \ \ \ \ \ \ \ \ \ \ \ \ \ \ \ \ \ \ \ \ \ \ \ \ $  
- Vitesse de convergence dans le t.l.c pour cha\^{\i}nes de Markov -}

\noindent{\large\bf I. INTRODUCTION}

\noindent {\it 
Dans ce papier on d\'esigne par $(E,\cE)$ un espace mesurable, par $Q$ une probabilit\'e de transition sur $(E,\cE)$ admettant 
une probabilit\'e invariante, not\'ee $\nu$, par $(X_n)_{n\geq0}$ une cha\^{\i}ne de Markov  sur $(E,\cE)$ associ\'ee \`a $Q$, et enfin 
par $\xi$ une fonction $\nu$-int\'egrable de $E$ dans $\R$ telle que $\nu(\xi) = 0$. 
On pose $\displaystyle S_n = \sum_{k=1}^{n} \xi(X_k)$. }

\noindent L'objet de ce travail est d'\'etudier la vitesse de convergence dans le th\'eor\`eme limite central (t.l.c) 
pour la suite de variables al\'eatoires $(\xi(X_n))_{n\geq0}$. La m\'ethode que nous utilisons s'inspire 
des techniques de perturbations d'op\'erateurs qui ont \'et\'e introduites par Nagaev \cite{nag1} \cite{nag2} 
et largement appliqu\'ees depuis. On trouvera dans \cite{hulo} un expos\'e g\'en\'eral de cette m\'ethode 
et de nombreuses r\'ef\'erences (voir \'egalement la remarque cons\'ecutive au th\'eor\`eme II).     

\noindent Les hypoth\`eses porteront sur le noyau $Q$ et les noyaux de Fourier $Q(t)$ associ\'es \`a $Q$ et $\xi$ ; 
en substance, on supposera que, sur un certain espace de Banach, $Q$ v\'erifie une hypoth\`ese de quasi-compacit\'e et 
que $Q(t)$ satisfait aux conditions du th\'eor\`eme de perturbations 
d'op\'erateurs de Keller-Liverani \cite{keli}. Comme il appara\^{\i}t d\'ej\`a dans \cite{ann} 
(dans le cadre des mod\`eles it\'eratifs) et dans \cite{ihp} (en vue d'obtenir un th\'eor\`eme local), 
le th\'eor\`eme de Keller-Liverani remplace avantageusement les \'enonc\'es classiques  
de perturbations d'op\'erateurs, notamment lorsque l'espace d'\'etats $E$ est non compact ou que $\xi$ est non born\'ee. 
Ainsi, dans \cite{ann} et \cite{ihp}, nous avons pu \'etablir des th\'eor\`emes limites sous des hypoth\`eses de moments polynomiaux 
l\`a o\`u habituellement \'etaient requis des moments exponentiels. 

\noindent Les r\'esultats de vitesse de convergence dans le t.l.c sont pr\'esent\'es au paragraphe II. 
On appliquera ensuite ($\S$ III) ces r\'esultats dans le cadre des cha\^{\i}nes 
$V$-g\'eom\'etriquement ergodiques et des mod\`eles it\'eratifs Lipschitziens. 

\noindent Pour illustrer les r\'esultats obtenus, consid\'erons l'exemple classique sur $\R^d$ 
des processus autor\'egressifs $X_n = A_nX_{n-1} + b_n$ sous les hypoth\`eses suivantes~: 
$A_1$ est presque-s\^urement une contraction stricte et $\E[\|b_1\|^3]<+\infty$. 
Alors, si $X_0$ a un moment d'ordre 2 et si 
$\xi$ est uniform\'ement lipschitzienne sur $\R^d$ (par exemple si $\xi(x)=\|x\|$), 
la vitesse de convergence dans le t.l.c. est en $n^{-\frac{1}{2}}$ (cf. Cor. III.2'). 
La condition de moment d'ordre 3 sur $b_1$ est la condition attendue en r\'ef\'erence au th\'eor\`eme de Berry-Esseen pour v.a.i.i.d. 

\noindent{\bf Notations.} Si $X$ et $Y$ sont des espaces de Banach, 
on d\'esigne par $X'$ le dual topologique de $X$, par $\cL(X)$ l'espace des endomorphismes continus de $X$,  
et par  $\cL(X,Y)$ l'espace des applications lin\'eaires continues de $X$ dans $Y$. Ces espaces 
sont munis des normes subordonn\'ees. On note $\langle\cdot,\cdot\rangle$ le crochet de dualit\'e sur $X'\times X$. \\ 
\noindent Pour $p\geq1$ on note $\L^p(\nu)$ 
l'espace de Lebesgue usuel associ\'e \`a $\nu$. On note ${\bf 1}=1_E$ la fonction identiquement \'egale \`a 1 sur $E$. \\
\noindent La probabilit\'e initiale de la cha\^{\i}ne sera appel\'ee $\mu_0$. 
La loi normale centr\'ee, de variance $\sigma^2$, est not\'ee ${\cal N}(0,\sigma^2)$. 
Enfin les noyaux de Fourier associ\'es \`a $Q$ et $\xi$ sont d\'efinis par 
$$ t\in\R,\ x\in E, \ \  Q(t)(x,dy) = e^{it\xi(y)}Q(x,dy).$$
\newpage

\noindent{\large\bf II. HYPOTH\`ESES ET \'ENONC\'ES DES R\'ESULTATS} 

\noindent 
Rappelons que $\nu$ d\'esigne une probabilit\'e $Q$-invariante. 
Les espaces de Banach sur lesquels on fera op\'erer $Q$ sont 
compos\'es de fonctions mesurables de $E$ dans $\C$. 
\'Etant donn\'e un tel espace $(\cB,\|\cdot\|)$, nous dirons que {\it $Q$ est $\cB$-g\'eom\'etriquement ergodique si~: 

\noindent ${\bf 1}\in\cB$, $\nu\in\cB'$, $Q\in\cL(\cB)$, et il existe $\kappa_0<1$, $C\geq0$ tels que }
$$\forall n\geq 1,\ \ \forall f\in \cB,\ \ \ \ \|Q^nf-\nu(f){\bf 1}\| \leq C\, \kappa_0^n\, \|f\|.$$
Sous cette condition, si $\xi\in\cB$, alors la s\'erie $\sum_{n\geq0} Q^n\xi$ converge 
absolument dans $\cB$ car $\nu(\xi)=0$. Dans ce cas on pose $\displaystyle \breve\xi = \sum_{n=0}^{+\infty} Q^n\xi$, et si en outre 
$\cB\subset \L^2(\nu)$, on note 
$$\sigma^2 = \nu(\breve\xi^2) - \nu((Q\breve\xi)^2),\ \ \mbox{et}\ \ 
\psi = Q(\breve\xi^2) - (Q\breve\xi)^2 - \sigma^2\, \bf{1}.$$
\noindent {\bf Hypoth\`ese ($\cH$).} {\it 
Il existe un espace de Banach $(\cB,\|\cdot\|)$ tel que $\cB\subset \L^3(\nu)$, $\xi\in\cB$, 
et v\'erifiant en outre les conditions suivantes~:  
 
\noindent {\bf (H1)} $Q$ est $\cB$-g\'eom\'etriquement ergodique. 

\noindent {\bf (H2)} On a $\displaystyle\sum_{p=0}^{+\infty} \nu(|Q^p\psi|^{\frac{3}{2}})^{\frac{2}{3}} < +\infty$. \\
\noindent {\bf (H3)} On a $\displaystyle \sup\bigg\{\nu(|e^{it\xi}-1|\, |f|), f\in\cB,\ \|f\|\leq 1\bigg\} = O(|t|)$ quand $t\r0$.  

\noindent {\bf (H4)}  Il existe un intervalle ouvert $I$ contenant $t=0$ tel que, pour $t\in I$, on ait $Q(t)\in\cL(\cB)$, et 
il existe des constantes $\kappa < 1,\ C\geq 0$ tels que 
$$\forall n\geq1,\ \ \forall t\in I,\ \ \forall f\in\cB,\ \ \ \ \|Q(t)^nf\| \leq C\, \kappa^n\, \|f\| + C\nu(|f|),$$
et enfin le rayon spectral essentiel de $Q(t)$ op\'erant sur $\cB$ 
est $\leq$ au r\'eel  $\kappa$. } 

\noindent 
On observera que la seule condition de moment sur $\xi$ est $\nu(|\xi|^3)<+\infty$. 

\noindent{\bf Th\'eor\`eme I.} {\it  
Supposons que $\mu_0 = \nu$ et $\sigma^2>0$. Sous l'hypoth\`ese (\cH), la suite $(\frac{S_n}{\sqrt n})_n$ 
converge en loi vers $\cN(0,\sigma^2)$ avec une vitesse de convergence au moins en $n^{-\frac{\tau}{2}}$ pour tout $\tau<1$, 
$$\mbox{\`a savoir~:}\ \ \ 
\forall \tau <1,\ \ \sup_{x\in\R}\bigg|\P(\frac{S_n}{\sigma\sqrt n}\leq x) - \cN(0,1)(]-\infty,x])\bigg|\ =\ 
O(n^{-\frac{\tau}{2}}).$$

\noindent {\bf Hypoth\`ese ($\widetilde{\cH}$).} Le couple $(Q,\cB)$ v\'erifie les hypoth\`eses ($\cH$), $\cB$ est contenu dans 
un espace de Banach $(\widetilde{\cB},\|\cdot\|_{\sim})$ qui s'envoie contin\^ument dans $\L^1(\nu)$, $Q$ est 
$\widetilde{\cB}$-g\'eom\'etriquement ergodique, et enfin 
$\displaystyle \sup\bigg\{\|Q(t)f-Qf\|_{\sim}, f\in\cB,\ \|f\|\leq 1\bigg\} = O(|t|)$ quand $t\r0$. 

\noindent{\bf Th\'eor\`eme II.}  
Supposons que l'hypoth\`ese ($\widetilde{\cH}$) soit satisfaite, et que $\sigma^2>0$ et 
$\mu_0\in\cB'\cap\widetilde{\cB}'$. Alors la vitesse de convergence dans le t.l.c est en $n^{-\frac{1}{2}}$. }

\noindent Sous la condition (H1), si $t\mapsto Q(t)$ est de classe $\cC^3$ de $I$ dans $\cL(\cB)$, 
le th\'eor\`eme de perturbations standard permet d'\'etablir une vitesse en $n^{-\frac{1}{2}}$ dans le t.l.c. \cite{gui}. 
La condition de r\'egularit\'e ci-dessus est par exemple v\'erifi\'ee lorque $\cB$ 
est une alg\`ebre de Banach (en g\'en\'eral dans ce cas $E$ est compact) et que $\xi\in\cB$ ; en effet  $t\mapsto Q(t)$ 
est alors analytique. 

\noindent {\bf Remarques.}\\
\noindent 1. La condition (H1) est \'equivalente au fait que $Q$ est quasi-compact sur $\cB$, avec 1 comme 
valeur propre simple et dominante. 

\noindent 2. Supposons que $Q$ soit $\cB$-g\'eom\'etriquement ergodique et que $\cB\subset \L^3(\nu)$, $\xi\in\cB$. 
Soit $(\cB_2,\|\cdot\|_{_2})$ un espace de Banach contenant les fonctions $g^2$, $g\in\cB$. 
S'il existe $A>0$ tel que l'on ait, pour $f\in\cB_2$, $\nu(|f|^{\frac{3}{2}})^{\frac{2}{3}} \leq A\|f\|_{_2}$, et si $Q$ est 
$\cB_2$-g\'eom\'etriquement ergodique, alors on a (H2). 
En effet $\breve\xi$, $Q\breve\xi$ et ${\bf 1}$ sont dans $\cB$, donc $\psi\in\cB_2$. Par 
d\'efinition du nombre $\sigma^2$ et par invariance de $\nu$, on a $\nu(\psi)=0$, par cons\'equent 
$\sum_{n\geq0} \|Q^n\psi\|_{_2} < +\infty$, et (H2) en d\'ecoule. 

\noindent 3. Soit  $(\cB_\gamma)_{0<\gamma\leq\gamma_0}$ une famille d'espaces de Banach, croissante pour l'inclusion. 
Les \'enonc\'es I-II sont particuli\`erement bien adapt\'es lorsque, pour tout $\gamma\in]0,\gamma_0]$, 
$Q$ est $\cB_\gamma$-g\'eom\'etriquement ergodique et que $Q(t)$ v\'erifie (H3)-(H4) sur $\cB_\gamma$ pour $|t|$ petit.  
Dans ce cas la principale difficult\'e r\'eside dans le choix de param\`etres 
$\gamma_1<\gamma_2<\gamma_3$ tels que, si $\cB=\cB_{\gamma_1}$, alors $Q$ v\'erifie sur 
$\cB_2 =\cB_{\gamma_2}$ les conditions de la remarque 2 (afin d'obtenir (H2)), et v\'erifie ($\widetilde{\cH}$) avec 
$\widetilde{\cB} = \cB_{\gamma_3}$. Les exemples du paragraphe III font intervenir de telles familles d'espaces. 

\noindent 4. Si, pour $t\in I$, l'in\'egalit\'e de (H4) est satisfaite et si l'ensemble $Q(t)(\{\|f\|\leq 1,\, f\in\cB\})$ est 
relativement compact dans $(\cB,\nu(|\cdot|))$, alors $Q(t)$ est quasi-compact \`a it\'er\'es born\'es sur $\cB$ \cite{itm}, et 
la propri\'et\'e dans (H4) sur le rayon spectral essentiel de $Q(t)$ est alors automatiquement satisfaite \cite{hen}. 
Cette remarque appliqu\'ee avec $t=0$ permet dans certain cas d'\'etablir (H1) lorsque $P$ v\'erifie en outre des hypoth\`eses 
d'irr\'eductibilit\'e et d'ap\'eriodicit\'e garantissant que 1 est une valeur propre 
simple et l'unique valeur propre de module 1 de $Q$. 

\noindent 5. Si $\widetilde{\cB}=\cB$, alors l'hypoth\`ese ($\widetilde{\cH}$) se r\'eduit \`a ($\cH$) et 
$\|Q(t) - Q\|_{\cL(\cB)} = O(|t|)$. Cependant cette derni\`ere condition est rarement satisfaite 
lorsque $\xi$ est non born\'ee, voir $\S$ II. 

\noindent 6. Une vitesse de convergence en $n^{-\frac{\tau}{4}}$ pour tout $\tau<1$ a \'et\'e 
\'etablie dans \cite{cr} pour des 
cha\^{\i}nes de Markov stationnaires, \`a espace d'\'etats compact, associ\'ees \`a un op\'erateur de transfert markovien. 
Dans \cite{cr} il n'est pas suppos\'e que $Q$ a une action quasi-compacte, et les conditions sur $\xi$ 
sont assez faibles. Mentionnons que la fonction $\psi$ introduite au d\'ebut du paragraphe 
est d\'ej\`a utilis\'ee dans les arguments de \cite{cr}. 

\noindent 7. Gr\^ace \`a (H1), la condition $\xi\in\cB$ est commode pour d\'efinir $\breve\xi$, et la condition 
$\cB\subset \L^3(\nu)$ assure alors que $\breve\xi\in\L^3(\nu)$. Cette derni\`ere propri\'et\'e n'est utilis\'ee 
qu'au $\S$ IV.1. \\
Si $\xi\notin\cB$, les th\'eor\`emes I-II subsistent, \`a condition de renforcer la condition (H2) comme suit~: 
la s\'erie $\sum_{n\geq0} Q^n\xi$ converge dans $\L^3(\nu)$, et la fonction $\psi$, que l'on peut alors d\'efinir comme au d\'ebut 
du paragraphe, v\'erifie $\sum_{p=0}^{+\infty} \nu(|Q^p\psi|^{\frac{3}{2}})^{\frac{2}{3}} < +\infty$.  

\noindent La preuve des th\'eor\`emes I- II est pr\'esent\'ee au $\S$ IV. 
La m\'ethode, fond\'ee sur des techniques de perturbations d'op\'erateurs introduites par Nagaev \cite{nag1} \cite{nag2},  
est proche de celle utilis\'ee dans \cite{ihp}. Plus pr\'ecis\'ement, on appliquera dans un 
premier temps une m\'ethode de r\'eduction en diff\'erence de martigale pour \'etablir que 
$\sup_{|t|\leq \sqrt n}|t|^{-1}\, |\E[e^{it\frac{S_n}{\sigma\sqrt n}}] - e^{\frac{-t^2}{2}}| = 
O(\frac{1}{\sqrt n})$. 
Cette majoration, obtenue en s'inspirant de \cite{jan}, permettra alors d'\'ecrire un d\'eveloppement limit\'e 
de la valeur propre dominante perturb\'ee de $Q(t)$ fournie par \cite{keli}. On appliquera alors 
les techniques usuelles de transform\'ee de Fourier. 

\noindent{\large\bf III. EXEMPLES} 

\noindent Dans les deux exemples trait\'es dans ce paragraphe, on d\'esigne par $(E,d)$ un espace m\'etrique 
non compact tel que toute boule ferm\'ee de $E$ soit compacte. On munit $E$ de sa tribu bor\'elienne $\cE$, et l'on note $x_0$ un point 
quelconque de $E$.  

\noindent{\bf III.1. Application aux cha\^{\i}nes $V$-g\'eom\'etriquement ergodiques.}

\noindent Soit $V$ une fonction mesurable de $E$ dans $[1,+\infty[$ telle que $V(x)\r+\infty$ quand $d(x,x_0)\r+\infty$, et 
soit $(X_n)_{n\geq0}$ une cha\^{\i}ne $V$-g\'eom\'etriquement ergodique 
\cite{mey} (chap. 16), \`a savoir~: il existe une probabilit\'e $Q$-invariante, $\nu$, telle que $\nu(V) < +\infty$, 
et $Q$ est $\cB_{_V}$-g\'eom\'etriquement ergodique, o\`u $\cB_{_V}$ est l'espace des fonctions mesurables sur $E$, \`a valeurs 
complexes, v\'erifiant 
$$\|f\|_{_{_V}} = \sup\{V(x)^{-1}|f(x)|,\ x\in E\} < +\infty.$$
Si $\xi^2$ est domin\'ee par un multiple de $V$, alors $(\frac{S_n}{\sqrt n})_n$ converge en loi vers $\cN(0,\sigma^2)$ \cite{mey}. 

\noindent{\bf Corollaire III.1.} {\it Si $|\xi|^3$ est domin\'ee par un multiple de $V$, 
si $\sigma^2>0$ et si $\E[V(X_0)] <+\infty$, alors 
$(\xi(X_n))_{n\geq0}$ v\'erifie le t.l.c avec une vitesse en $n^{-\frac{1}{2}}$. }

\noindent Dans le cas stationnaire et pour $V={\bf 1}$, nous retrouvons la vitesse en $n^{-\frac{1}{2}}$ fournie par 
le th\'eor\`eme de Bolthausen \cite{bolt} (les conditions de \cite{bolt} sur les coefficients de m\'elange sont clairement 
satisfaites lorsque $V={\bf 1}$). Dans le cadre sp\'ecifique des cha\^{\i}nes $V$-g\'eom\'etriquement ergodiques, 
une vitesse en $\frac{1}{\sqrt n}$ est obtenue  dans \cite{konmey} sous la condition assez restrictive que $\xi$ est born\'ee. 
Dans \cite{fuh} la vitesse, exprim\'ee en termes d'in\'egalit\'es de 
Paley, est pr\'esent\'ee dans le cas stationnaire sous des conditions de moment assez fortes.  
Enfin \cite{stein} \'etablit une vitesse en $(\frac{\ln n}{n})^{\beta}$ lorsque $\xi$ est domin\'ee par un multiple de $V^\alpha$ 
($0<\alpha\leq\frac{1}{2}$), avec $\beta = \frac{1}{2(\alpha+1)}$. \\ 
\noindent Ainsi le corollaire III.1 am\'eliore les r\'esultats des travaux pr\'ec\'edemment 
cit\'es. En outre, comme $\nu(V) < +\infty$, 
la condition que $|\xi|^3$ est domin\'ee par un multiple de $V$ est proche de celle requise dans le th\'eor\`eme de Berry-Esseen 
puisqu'elle exprime d'une certaine fa\c con que $\xi$ admet un moment d'ordre 3. 
\noindent Mentionnons enfin que, gr\^ace \`a la m\'ethode spectrale et au th\'eor\`eme de Keller-Liverani, 
des th\'eor\`emes limites local et de renouvellement ont \'et\'e \'etablis dans \cite{prepub} 
lorsque $\xi$ est non-arithm\'etique et domin\'ee par un multiple de $V^{\frac{1}{2}-\varepsilon}$ ($\varepsilon>0$).  

\noindent {\it Preuve du corollaire III.1.} 
Soit $W=V^{\frac{1}{3}}$, $U=V^{\frac{2}{3}}$, et soient $\cB_{_W}$, $\cB_{_U}$ les espaces obtenus 
comme ci-dessus en rempla\c cant $V$ respectivement par $W$ et $U$. Notons que $\nu\in\cB_{_W}'$ car $\nu(W)<+\infty$, que 
$\xi\in\cB_{_W}$, et enfin que $\cB_{_W}$ s'envoie contin\^ument dans $\L^3(\nu)$ car $\nu(W^3) = \nu(V)<+\infty$. 
Nous allons montrer que $Q$ v\'erifie l'hypoth\`ese ($\widetilde{\cH}$) avec $\cB=\cB_{_W}$, $\widetilde{\cB}=\cB_{_V}$. \\
Comme $Q$ est $V$-g\'eom\'etriquement ergodique par hypoth\`ese, il l'est \'egalement relativement \`a $W$,  
voir \cite{mey} (Lem. 15.2.9 et Th. 16.0.1). D'o\`u (H1). Pour les m\^emes raisons $Q$ est $U$-g\'eom\'etriquement ergodique. 
En outre $\cB_{_U}$ s'envoie contin\^ument dans $L^{\frac{3}{2}}(\nu)$ car $\nu(U^{\frac{3}{2}}) = \nu(V)<+\infty$. 
D'o\`u (H2) (cf. Rq. 2 avec $\cB_2 =\cB_{_U}$). \\
On a clairement $Q(t)\in\cL(\cB_{_W})$ pour tout $t\in\R$. Comme $\cB_{_W}$ s'envoie contin\^ument dans $\L^3(\nu)$, 
l'in\'egalit\'e de H\"older implique (H3). 
La propri\'et\'e (H4), plus difficile \`a \'etablir, r\'esulte d'un travail r\'ecent de Hennion \cite{hub}, 
voir \cite{prepub}. Enfin si $f\in\cB_{_W}$, alors 
$$|Q(t)f-Qf| \leq Q\bigg(|e^{it\xi}-1|\, |f|\bigg) \leq |t|\ Q\left(|\xi|W\right)\, 
\|f\|_{_W} \leq |t|\, \|\xi\|_{_W}\, Q(W^2)\, \|f\|_{_W},$$
d'o\`u $\|Q(t)f-Qf\|_{_V}  \leq |t|\, \|\xi\|_{_W}\, \|QW^2\|_{_V}\, \|f\|_{_W}$. Ce qui pr\'ec\`ede prouve ($\widetilde{\cH}$). 
Comme par hypoth\`ese $\mu_0(V)<+\infty$, on a $\mu_0\in\cB_{_W}'\cap\cB_{_V}'$, et le corollaire d\'ecoule du th\'eor\`eme II. \fdem

\noindent{\bf III.2. Application aux Mod\`eles it\'eratifs Lipschitziens.} 

\noindent On d\'esigne par $G$  un semi-groupe de transformations lipschitziennes de $E$, et 
par  $\cal G$ une tribu sur $G$. On suppose que l'action de $G$ sur $E$ est mesurable. 
Pour $g\in G$, on pose 
$$c(g) = \sup\bigg\{\frac{d(gx,gy)}{d(x,y)},\ x,y\in E,\ x\neq y\bigg\}.$$
\noindent Soit $\displaystyle (Y_n)_{n\geq 1}$ une suite de v.a.i.i.d. \`a valeurs dans $G$. On note $\pi$ leur loi commune. 
\'Etant donn\'ee une variable al\'eatoire $X_0$ \`a valeurs dans $E$, ind\'ependante de $(Y_n)_n$, de loi $\mu_0$, 
on consid\`ere la suite 
$(X_n)_{n\geq0}$ d\'efinie pour $n\geq 1$ par $X_n = Y_nX_{n-1}$. 

\noindent Alors $(X_n)_{n\geq0}$ est une cha\^{\i}ne de Markov de probabilit\'e de transition 
$\displaystyle (Qf)(x) = \int_G f(gx)d\pi(g)$. 

\noindent On suppose dans ce paragraphe qu'il existe une constante $C\geq 0$ telle que l'on ait 
$$\forall (x,y)\in E^2,\ \ \ \ |\xi(x)-\xi(y)| \leq C\, d(x,y).$$
\noindent Sous les hypoth\`eses $\int_G c(g)^2 d\pi(g) < 1$ et $\int_G d(gx_0,x_0) ^2 d\pi(g) +\infty$, 
la suite $(\frac{S_n}{\sqrt n})_n$ converge en loi vers une gaussienne ${\cal N}(0,\sigma^2)$, voir \cite{ben}. 

\noindent Soit $\Gamma(g) = 1 + c(g) + d(g x_0,x_0)$. On renforce ici la condition pr\'ec\'edente en supposant qu'il existe un entier 
$n_0\geq1$ tel que  
$${\bf (*)} \ \ \ \ \ \ \ \ \ \int_G \Gamma(g)^3\,  (1+c(g)^{\frac{1}{2}})\, d\pi(g) < +\infty\ \ \ \mbox{et} 
\ \ \ \ \int_G c(g)^{\frac{1}{2}}\, \max\{c(g),1\}^3\, d\pi^{*n_0}(g) < 1,$$
o\`u l'on a d\'esign\'e par $\pi^{*n_0}$ la loi de $Y_{n_0}\cdots Y_1$. \\
Sous ces hypoth\`eses, on sait qu'il existe une unique probabilit\'e $Q$-invariante, $\nu$, 
et que $\nu(d(\cdot,x_0)^{3})<+\infty$ (appliquer le Th. I de \cite{ann} avec la distance $d(x,y)^\frac{1}{2}$).  

\noindent{\bf Corollaire III.2.} {\it 
Sous l'hypoth\`ese $(*)$, si $ \sigma^2>0$ et $\E[d(X_0,x_0)^{2}] < +\infty$, alors $(\xi(X_n))_{n\geq0}$ v\'erifie le t.l.c 
avec une vitesse en $n^{-\frac{1}{2}}$. }

\noindent Le cadre ci-dessus contient celui des mod\`eles it\'eratifs \cite{duf}, en particulier celui des processus autor\'egressifs 
d\'efinis par une v.a $X_0$ \`a valeurs dans $\R^d$, puis par $X_n = A_nX_{n-1} + b_n$ ($n\geq1$), o\`u 
$(A_n,b_n)_{n\geq1}$ est une suite de v.a.i.i.d \`a valeurs dans $\cM_d(\R)\times\R^d$, ind\'ependante de $X_0$ ; on a not\'e 
$\cM_d(\R)$ l'espace des matrices r\'eelles carr\'ees d'ordre $d$. \\
Dans ce contexte une vitesse en $\frac{1}{\sqrt n}$ dans le t.l.c a \'et\'e \'etablie dans \cite{mira} lorsque $A_1$ est 
presque-s\^urement une contraction stricte et que $b_1$ v\'erifie une condition de moment exponentiel. Dans \cite{cun} 
une vitesse en $n^{-\frac{\tau}{2}}$ ($\tau<1$) est obtenue pour une large classe de fonctions $\xi$, sous la m\^eme 
condition sur $A_1$ et sous l'hypoth\`ese $\E[\|b_1\|^p]<+\infty$ pour tout $p\in\N$. Les r\'esultats de 
\cite{mira} ont \'et\'e \'etendus dans \cite{ann} sous une condition de contraction en moyenne sur $\|A_1\|$  et 
sous la condition de moment $\E[\|b_1\|^{8+\varepsilon}]<+\infty$ ($\varepsilon>0$), 
o\`u $\|\cdot\|$ d\'esigne indiff\'eremment une norme de $\R^d$ et la norme subordonn\'ee associ\'ee sur $\cM_d(\R)$. 
Du corollaire III.2 nous d\'eduisons par exemple le r\'esultat suivant. 

\noindent{\bf Corollaire III.2'.} {\it Si $\|A_1\| < 1$ presque s\^urement, si $\E[\|b_1\|^3]<+\infty$ et 
$\E[\|X_0\|^2] < +\infty$, alors $(\xi(X_n))_{n\geq0}$ v\'erifie le t.l.c avec une vitesse en $n^{-\frac{1}{2}}$ 
(sous r\'eserve que $\sigma^2>0$). }

\noindent{\it Preuve du corollaire III.2.} 
Soit $\lambda\in]0,1]$ quelconque, soit $p_\lambda(x) = 1+\lambda\, d(x,x_0)^\frac{1}{2}$ (le r\'eel 
$\lambda$ sera choisi ult\'erieurement), 
et pour $\gamma>0$, soit $\cL_\gamma$ l'espace de Banach 
des fonctions $f$ de $E$ dans $\C$telles que 
$$m_\gamma(f) = \sup\bigg\{\frac{|f(x)-f(y)|}{d(x,y)^\frac{1}{2} 
p_\lambda(x)^\gamma  p_\lambda(y)^\gamma},\ x,y\in E,\ x\neq y\bigg\}\ < +\infty,$$
muni de la norme $\|f\|_\gamma = m_\gamma(f) + |f|_\gamma $, o\`u l'on a pos\'e 
$\displaystyle \ |f|_\gamma  = \sup_{x\in E}\ \frac{|f(x)|}{p_\lambda(x)^{\gamma+1}} < +\infty$. \\ 
On a clairement $\xi\in\cL_1$, et si $f\in\cL_1$, alors $\nu(|f|^3)^{\frac{1}{3}} \leq |f|_1\, 
\nu(p_\lambda^{6})^{\frac{1}{3}}$, avec $\nu(p_\lambda^{6})<+\infty$. Donc $\cL_1$ s'envoie contin\^ument dans $\L^3(\nu)$.

\noindent{\bf Lemme III.1. } {\it On a (H1) avec $\cB=\cL_1$, et (H2). }

\noindent{\it Preuve.} On a $\nu\in\cL_1'\cap\cL_3'$ car $\nu(p_\lambda^{4})<+\infty$.  
En choisissant $\lambda$ suffisamment petit, l'ergodicit\'e g\'eom\'etrique de $Q$ relativement \`a $\cB=\cL_1$, puis \`a 
$\cL_3$, r\'esulte des conditions ($*$) et de \cite{ann} (Th. 5.5 appliqu\'e avec la distance $d(x,y)^\frac{1}{2}$) 
\footnote{\`A cet effet le r\'eel $\lambda$ doit \^etre fix\'e tel que $\pi(c^\frac{1}{2}\delta_\lambda^6) < 1$, o\`u  
$\delta_\lambda(g) = \max\{c(g),1\}^\frac{1}{2} + \lambda\, d(gx_0,x_0)^\frac{1}{2}$, ce qui est possible gr\^ace aux conditions $(*$) 
et au th\'eor\`eme de convergence domin\'ee.}.   
En outre si $f\in\cL_3$, alors $\nu(|f|^{\frac{3}{2}})^{\frac{2}{3}} \leq |f|_3\, \nu(p_\lambda^{6})^{\frac{2}{3}} \leq 
\|f\|_3\, \nu(p_\lambda^{6})^{\frac{2}{3}}$, et $\cL_3$ contient les fonctions $g^2$, $g\in\cL_1$. D'o\`u (H2) (cf. Rq. 2 
avec $\cB_2 = \cL_3$). \fdem 
 
\noindent{\bf Lemme III.2.} {\it La condition (H3) est satisfaite, $Q(t)$ est un endomorphisme continu de $\cL_1$ pour tout $t\in\R$, 
et enfin on a (H4) pour $|t|$ petit. }

\noindent{\it Preuve.} 
Comme $\cL_1$ s'envoie contin\^ument dans $\L^3(\nu)$, (H3) r\'esulte de l'in\'egalit\'e de H\"older. \\
On pose $\delta_\lambda(g) = \max\{c(g),1\}^\frac{1}{2} + \lambda\, d(gx_0,x_0)^\frac{1}{2}$. 
On a $\displaystyle \sup_{x\in E} \frac{p_\lambda(gx)}{p_\lambda(x)} \leq \delta_\lambda(g)$ et 
$\delta_\lambda(g)^2 \leq 2\Gamma(g)$. Soit $f\in\cL_1$. De la d\'efinition 
de $c(g)$ et de l'in\'egalit\'e $|e^{ia}-e^{ib}| \leq 2|b-a|^\frac{1}{2}$ ($a,b\in\R$), on obtient en posant 
$A=\pi(c^\frac{1}{2}\delta_\lambda^2)$ 
\begin{eqnarray*}
|(Q(t)f)(x) - (Q(t)f)(y)| &\leq& \int |f(gx)-f(gy)|\, d\pi(g)\ +\ 
\int |f(gy)|\, |e^{it\xi(gx)} - e^{it\xi(gy)}| \, d\pi(g) \\
&\leq& A \, m_1(f) d(x,y)^\frac{1}{2} p_\lambda(x)  p_\lambda(y) \ +\ 
2AC^\frac{1}{2}|t|^\frac{1}{2} |f|_1  d(x,y)^\frac{1}{2} p_\lambda(y)^2. 
\end{eqnarray*}
On peut \'evidemment supposer que $d(y,x_0) \leq d(x,x_0)$ (sinon, inverser le r\^ole jou\'e par $x$ et $y$), de sorte que 
$p_\lambda(y)^2 \leq p_\lambda(x) p_\lambda(y)$, et il vient que $Q(t)f\in\cL_1$, avec 
$$m_1(Q(t)f) \leq  A\, m_1(f) \ +\ 2AC^\frac{1}{2}\, |t|^\frac{1}{2} |f|_1.$$
On a clairement $|Q(t)f|_1 \leq |f|_1\, \pi(\delta_\lambda^2)$, 
de sorte que $Q(t)$ a une action continue sur $\cL_1$. \\ 
Pour \'etablir l'in\'egalit\'e de (H4), on utilise le fait que les normes $\|\cdot\|_1$ et 
$\|\cdot\|_\nu = m_1(\cdot)+\nu(|\cdot|)$ 
sont \'equivalentes \cite{ann} (Sect. 5). Alors, d'apr\`es l'in\'egalit\'e ci-dessus, il existe une constante $D>0$ telle que 
$$m_1(Q(t)f) \leq  A\, m_1(f) \ +\ D|t|^\frac{1}{2}\, [m_1(f)+\nu(|f|)] =   
(A+ D|t|^\frac{1}{2}) \, m_1(f) \ +\ D|t|^\frac{1}{2}\, \nu(|f|).$$
On a $\nu(|Q(t)f|) \leq \nu(Q|f|) = \nu(|f|)$. Donc $\|Q(t)f\|_\nu \leq (A+ D|t|^\frac{1}{2}) \|f\|_\nu +  
(D|t|^\frac{1}{2}+1)\, \nu(|f|)$. Le r\'eel $\lambda$ fix\'e dans la preuve du lemme III.1 est tel que 
$A = \pi(c^\frac{1}{2}\delta_\lambda^2) < 1$. 
Soit alors $t_0$ tel que $\kappa = A+ D|t_0|^\frac{1}{2} < 1$, soit $C' = D|t_0|^\frac{1}{2}+1$, et soit $t$ tel que $|t|\leq t_0$. 
Comme $\nu(|Q(t)f|) \leq \nu(Q|f|) = \nu(|f|)$, il r\'esulte d'une r\'ecurrence \'evidente que 
$\|Q(t)^nf\|_\nu \leq \kappa^n\, \|f\|_\nu + \frac{C'}{1-\kappa}\nu(|f|)$, ce qui prouve l'in\'egalit\'e de (H4). 

\noindent Il reste \`a \'etablir dans (H4) la propri\'et\'e relative au rayon spectral essentiel de $Q(t)$. 
On a $\nu(|Q(t)^nf|) \leq \nu(Q^n|f|) = \nu(|f|)$, et la boule unit\'e de $(\cL_1,\|\cdot\|_\nu)$ est relativement 
compacte dans  $(\cL_1,\nu(|\cdot|))$ (utiliser le th\'eor\`eme d'Ascoli et le th\'eor\`eme de Lebesgue). 
La  propri\'et\'e souhait\'ee r\'esulte alors de l'in\'egalit\'e de (H4) et de  \cite{hen}. \fdem 

\noindent Ce qui pr\'ec\`ede montre que $(Q,\cL_1)$ v\'erifie ($\cH$). 
D\'emontrons maintenant que $Q$ v\'erifie ($\widetilde{\cH}$) avec $\widetilde{\cB} = \cL_3$. 
On a d\'ej\`a vu que $Q$ est $\cL_3$-g\'eom\'etriquement ergodique. En outre on a~:  

\noindent{\bf Lemme III.3.} {\it Il existe $E>0$ tel que, pour $f\in\cL_1$, $t\in\R$, on ait 
$\|Q(t)f - Qf\|_3 \leq E|t|\|f\|_1$.} 

\noindent{\it Preuve.} Soit $f\in\cL_1$. Il existe $D>0$ tel que $\xi \leq D\, p_\lambda^2$. D'o\`u  
$$|(Q(t)f)(x) - Qf(x)| \leq \int |e^{it\xi(gx)}-1|\, |f(gx)|\, d\pi(g)  
\leq D\, |f|_1\, |t|\, p_\lambda(x)^4\, \pi(\delta_\lambda^4),$$
donc $|Q(t)f - Qf|_3 \leq  D\, |f|_1\, |t|\, \pi(\delta_\lambda^4)$. En outre, en posant 
$\tilde A=\pi(c^\frac{1}{2}\delta_\lambda^4)$, $\tilde B= \pi(c \delta_\lambda^2)$, on a 
\begin{eqnarray*}
\bigg|[(Q(t)f)(x) - Qf(x)] &-& [(Q(t)f)(y)-Qf(y)]\bigg| \\
&\leq& \tilde A\, D\, m_1(f) |t| d(x,y)^\frac{1}{2} p_\lambda(x)^3  p_\lambda(y) \ +\ 
\tilde B C|t|\, |f|_1  d(x,y) p_\lambda(y)^2.  
\end{eqnarray*}
On a $p_\lambda(y) \leq p_\lambda(y)^3$ et $d(x,y)^\frac{1}{2} \leq\frac{1}{\lambda}(p_\lambda(x)+p_\lambda(y)) \leq 
\frac{2}{\lambda} p_\lambda(x) p_\lambda(y) \leq \frac{2}{\lambda} p_\lambda(x)^3 p_\lambda(y)$, 
par cons\'equent $m_3(Q(t)f - Qf) \leq (\tilde A\, D\, m_1(f) + \frac{2}{\lambda}\tilde B C |f|_1)\, |t|$. \fdem 

\noindent Comme $\mu_0(p_\lambda^4) < +\infty$ par hypoth\`ese, on a $\mu_0\in\cL_1'\cap\cL_3'$, et le corollaire III.2 r\'esulte 
du th\'eor\`eme II. \fdem 

\newpage

\noindent{\large\bf IV. D\'EMONSTRATION DU TH\'EOR\`EME} 

\noindent{\bf IV.1. Une in\'egalit\'e sur les fonctions caract\'eristiques}

\noindent On suppose dans ce paragraphe que $X_0$ suit la loi $\nu$, que les conditions (H1)-(H2) sont satisfaites, 
et enfin que $\cB\subset \L^3(\nu)$, $\xi\in\cB$. 
Les \'el\'ements $\breve\xi$, $\sigma^2$ et $\psi$ ont \'et\'e d\'efinis au d\'ebut du paragraphe II, et pour $n\geq 1$ on pose 
$$U_n = \breve\xi(X_{n}) - Q\breve\xi(X_{n-1})\ \ \ \mbox{et}\ \ \ \ T_n = U_1+\cdots+U_n.$$
\noindent {\bf Proposition IV.1.} {\it Si $\sigma^2>0$, alors il existe une constante $C>0$ telle que l'on ait }
$$\forall n\in\N^*,\ \ \forall t\in[-\sqrt n,\sqrt n],\ \ \ \ \ 
\bigg|\E[e^{it\frac{T_n}{\sigma\sqrt n}}]\, - \, e^{-\frac{t^2}{2}}\bigg| \leq C\frac{|t|}{\sqrt n}.$$
Admettons pour le moment cette proposition. En utilisant la m\'ethode de Gordin \cite{gor}, nous allons 
en d\'eduire le r\'esultat suivant. 

\noindent {\bf Corollaire IV.1.} {\it  Si $\sigma^2>0$, alors il existe une constante $C>0$ telle que l'on ait }
$$\forall n\in\N^*,\ \ \forall t\in[-\sqrt n,\sqrt n],\ \ \ \ \
\bigg|\E[e^{it\frac{S_n}{\sigma\sqrt n}}]\, - \, e^{-\frac{t^2}{2}}\bigg| \leq 
C\frac{|t|}{\sqrt n}.$$ 
\noindent{\it Preuve du corollaire.} Soit $V_n= Q\breve\xi(X_0) - Q\breve\xi(X_{n})$. Gr\^ace \`a 
l'\'equation de Poisson $\breve\xi -Q\breve\xi = \xi$, on obtient que $S_n = T_n + V_n$. Par ailleurs, 
de la stationnarit\'e de $(X_n)_n$ et du fait que $\breve\xi\in\cB\subset \L^1(\nu)$, il vient que $\sup_n\E[|V_n|] < +\infty$. 
Enfin on a 
$$\bigg|\E[e^{it\frac{S_n}{\sigma\sqrt n}}]\, - \, e^{-\frac{t^2}{2}}\bigg| = 
\bigg|\E[e^{it\frac{T_n}{\sigma\sqrt n}}\, e^{it\frac{V_n}{\sigma\sqrt n}}]\, - \, e^{-\frac{t^2}{2}}\bigg| \leq 
\bigg|\E[e^{it\frac{T_n}{\sigma\sqrt n}}]\, - \, e^{-\frac{t^2}{2}}\bigg| + \E[|e^{it\frac{V_n}{\sigma\sqrt n}}-1|],$$ 
avec $\E[|e^{it\frac{V_n}{\sigma\sqrt n}}-1|] \leq \frac{1}{\sigma}\, \frac{|t|}{\sqrt n}\, \sup_n\E[|V_n|]$. 
On conclut alors gr\^ace \`a la proposition. \fdem

\noindent{\it Preuve de la proposition.} On note $\cF_n = \sigma(X_0,\ldots,X_n)$ pour $n\geq0$. Il est facile de voir que 
$(U_n)_{n\geq1}$ est une suite stationnaire d'accroissements de martingale relativement \`a $(\cF_n)_n$, telle que 
$\E[|U_1|^3]< +\infty$ car $\breve\xi\in\cB\subset \L^3(\nu)$. \\ 
Pour simplifier nous consid\'erons le cas $\sigma^2=1$, et nous adaptons au cadre markovien les majorations pr\'esent\'ees dans 
la preuve du th\'eor\`eme 6 de \cite{jan} [pp. 41-44]. \`A cet effet on pose $T_0=0$, $W_n = U_n^2-1$ pour $n\geq1$, et l'on rappelle 
que $e^{ix} = 1+ix-\frac{x^2}{2}+u(ix)$, avec $|u(ix)|\leq \frac{|x|^3}{6}$. \\ 
En \'ecrivant $\E[e^{i\lambda T_n}] = \E[e^{i\lambda T_{n-1}}\, e^{i\lambda U_n}]$, 
puis en appliquant la remarque pr\'ec\'edente avec $x=\lambda U_n$, et enfin en observant que 
$\E[e^{i\lambda T_{n-1}}\, U_n] = \E[e^{i\lambda T_{n-1}}\, \E[U_n \, |\, \cF_{n-1}]] = 0$, il est facile de voir par r\'ecurrence que 
pour $t\in\R$ et $n\geq1$
$$\E[e^{i\frac{t}{\sqrt n} T_n}] - e^{-\frac{t^2}{2}} = A_n(t) + B_n(t) + C_n(t)\ \ \ \ \ \mbox{avec}$$
$$A_n(t) = (1-\frac{t^2}{2n})^n - e^{-\frac{t^2}{2}},\ \ \ \mbox{puis} \ \ \ \ \ 
B_n(t) = \sum_{k=0}^{n-1}(1-\frac{t^2}{2n})^k \E\bigg[e^{i\frac{t}{\sqrt n}T_{n-k-1}}u(i\frac{t}{\sqrt n}U_{n-k})\bigg],$$
$$\mbox{et enfin}\ \ \ \ C_n(t) = -\frac{t^2}{2n} \sum_{k=0}^{n-1}(1-\frac{t^2}{2n})^k 
\E\bigg[e^{i\frac{t}{\sqrt n}T_{n-k-1}}\, W_{n-k}\bigg].$$ 
Soit $t\in[-\sqrt n,\sqrt n]$. On a 
$$-A_n(t) \leq  e^{-\frac{t^2}{2}} - e^{n(-\frac{t^2}{2n} -\frac{3t^4}{8n^2})} 
\leq (1-e^{-\frac{3t^4}{8n})})\, e^{-\frac{t^2}{2}} \leq  \frac{3t^4}{8n}\, e^{-\frac{t^2}{2}} \leq C_1\frac{|t|}{\sqrt n},$$
$$|B_n(t)| \leq \sum_{k=0}^{n-1}(1-\frac{t^2}{2n})^k\, \frac{|t|^3}{6n\sqrt n}\, \E[|U_{n-k}|^3] = 
\frac{1}{3}\E[|U_1|^3]\, \frac{|t|}{\sqrt n}.$$
Pour la majoration de $C_n(t)$, on utilise le lemme suivant 

\noindent {\bf Lemme IV.1.} {\it Pour $k\geq\ell\geq1$ on a $\E[W_k\, |\, \cF_{\ell-1}] = (Q^{k-\ell}\psi)(X_{\ell-1})$. }

\noindent{\it Preuve.} On a $W_k  = [\breve\xi(X_{k}) - Q\breve\xi(X_{k-1})]^2 - 1$. 
Le lemme s'\'etablit alors ais\'ement en d\'eveloppant cette expression et en utilisant 
le fait que $(X_n)_{n\geq0}$ est une cha\^{\i}ne de Markov. \fdem 

\noindent D'apr\`es (H2) on sait que $C_3 = \sum_{p=0}^{+\infty} \nu(|Q^p\psi|^{\frac{3}{2}})^{\frac{2}{3}} < +\infty$, donc 
$\sum_{p=0}^{+\infty} \nu(|Q^p\psi|) < +\infty$. Effectuons maintenant sur $W_\ell$ 
une r\'eduction en diff\'erence de martingale, \`a savoir $W_\ell = Y_\ell +  Z_\ell$ pour $\ell\geq 1$, avec 
$$Y_\ell = \sum_{p=0}^{+\infty}\bigg[\E[W_{p+\ell}\, |\, \cF_\ell] - \E[W_{p+\ell}\, |\, \cF_{\ell-1}]\bigg]\ \ \ \mbox{et}\ \ \ 
Z_\ell = \sum_{p=0}^{+\infty}\E[W_{p+\ell}\, |\, \cF_{\ell-1}] - \sum_{p=1}^{+\infty}\E[W_{p+\ell}\, |\, \cF_{\ell}].$$
$Y_\ell$ est $\cF_\ell$-mesurable, et $\E[Y_\ell\, |\, \cF_{\ell-1}] = 0$. Comme $T_{\ell-1}$ est $\cF_{\ell-1}$-mesurable, on 
obtient $\E[e^{itT_{\ell-1}}\, Y_\ell ] = \E[e^{itT_{\ell-1}}\,\E[Y_\ell |\, \cF_{\ell-1}]] = 0$, 
d'o\`u $\E[e^{itT_{\ell-1}}\, W_\ell ] = \E[e^{itT_{\ell-1}}\, Z_\ell ]$. Par cons\'equent, en posant 
$\displaystyle Z'_\ell = \sum_{p=0}^{+\infty}\E[W_{p+\ell}\, |\, \cF_{\ell-1}]$, il vient 
$$(E)\ \ \ \ \ \ \E[e^{itT_{\ell-1}}\, W_\ell ]\ =\ \E[e^{itT_{\ell-1}}\, Z'_\ell ] - \E[e^{itT_{\ell-1}}\, Z'_{\ell+1}].$$ 
On a $Z'_\ell = \sum_{p=0}^{+\infty} (Q^p\psi)(X_{\ell-1})$ (Lemme IV.1), et comme $\sum_{p=0}^{+\infty} \nu(|Q^p\psi|) < +\infty$, 
il vient 
\begin{eqnarray*}
\E[e^{itT_{\ell-1}}\, Z'_\ell] &=& \sum_{p=0}^{+\infty}\E\bigg[e^{itT_{\ell-1}}\, (Q^{p}\psi)(X_{\ell-1})\bigg] \\
\mbox{(par stationnarit\'e de}\ (X_n)_{n\geq0}) \ &=& \sum_{p=0}^{+\infty}\E\bigg[e^{it(U_2+\cdots+U_\ell)}\, (Q^{p}\psi)(X_\ell)\bigg]\\
\mbox{(par le lemme IV.1.)}\ &=& \sum_{p=0}^{+\infty}\E\bigg[e^{it(T_\ell-U_1)}\,  \E[W_{p+\ell+1}\, |\, \cF_{\ell}]\bigg] 
= \E[e^{it(T_\ell-U_1)}\, Z'_{\ell+1}].
\end{eqnarray*}
De l'\'egalit\'e (E) on d\'eduit que 
$$|\, \E[e^{itT_{\ell-1}}\, W_\ell ]\, | \leq \E[\, |Z'_{\ell+1}|\, |e^{it(U_\ell-U_1)}-1|\, ] 
\leq |t|\,  \E[|Z'_{\ell+1}|^{\frac{3}{2}}]^{\frac{2}{3}}\, \E[|U_\ell-U_1|^3]^{\frac{1}{3}}$$ 
avec $\E[|Z'_{\ell}|^{\frac{3}{2}}]^{\frac{2}{3}} \leq C_3$, puis par stationnarit\'e 
$\E[|U_\ell-U_1|^3]^{\frac{1}{3}} \leq 2\, \E[|U_1|^3]^{\frac{1}{3}}$. 
Par cons\'equent 
$|\E[e^{i\frac{t}{\sqrt n} T_{\ell-1}}\, W_\ell ]| \leq 2\, \E[|U_1|^3]^{\frac{1}{3}}\, C_3\, \frac{|t|}{\sqrt n}$, et finalement 
$\displaystyle |C_n(t)| \leq 2\, C_3\, \E[|U_1|^3]^{\frac{1}{3}}\, \frac{|t|}{\sqrt n}.$ \fdem 

\noindent{\bf IV.2. Un th\'eor\`eme de perturbations}

\noindent {\bf Th\'eor\`eme IV.} {\it Supposons que les conditions (H1) (H3) (H4)  soient satisfaites. \\
Soit $0<\tau<1$. Il existe un intervalle ouvert $J\subset I$ centr\'e en $t=0$, et des applications 
$\lambda(\cdot),\, v(\cdot),\, \phi(\cdot)$ et $N(\cdot)$ 
\`a valeurs respectivement dans $\C,\, \cB,\, \cB'$ et $\cL(\cB)$ tels que l'on ait, pour $t\in J$, $n\geq1$, et $f\in\cB$, 
$$(D)\ \ \ \ \ \ \ \ \ \ Q(t)^nf = \lambda(t)^n \langle \phi(t),f \rangle v(t) + N(t)^n f,$$ 
avec en outre les propri\'et\'es suivantes~: 

\noindent {\bf (a)} $\langle \phi(t),v(t)\rangle =1$, $\phi(t) N(t) = 0$, $N(t)v(t) = 0$,  $\displaystyle\lim_{t\r0}\lambda(t) = 1$. 

\noindent {\bf (b)} $\ \langle \nu,v(t) \rangle = 1$, et il existe $C>0$ et $\rho<1$ tels que l'on ait pour $t\in J$ et $n\geq 1$ \\
\indent {\bf (b1)} $\langle \nu,|v(t)-{\bf 1}| \rangle \leq C\, |t|^\tau$ \\
\indent {\bf (b2)} $|\langle \phi(t),{\bf 1} \rangle - 1 | \leq  C\, |t|^\tau$ \\
\indent {\bf (b3)} $\langle \nu, |N(t)^n {\bf 1}| \rangle \leq  C\, \rho^n\,|t|^\tau $. 

\noindent {\bf (c)} $\lambda(\cdot)$ est continue sur $J$. }

\noindent{\it Preuve.} Soit $t\in I$, $n\geq 1$, $f\in\cB$. On a $|Q(t)^nf| \leq Q^n|f|$, donc par invariance de $\nu$, 
$\nu(|Q(t)^nf|) \leq \nu(Q^n|f|) =\nu(|f|)$, puis pour $t_0\in I$ et $h\in \R$ tel que $t_0+h\in I$ 
$$\nu(|Q(t_0+h)f-Q(t_0)f|) \leq \nu(Q(|e^{ih\xi}-1|\, |f|) = \nu(|e^{ih\xi}-1|\, |f|).$$ 
D'o\`u $\sup\{\nu(|Q(t_0+h)f-Q(t_0)f|), f\in\cB,\, \|f\|\leq 1\} = O(|h|)$ d'apr\`es (H3). \\
La condition (H4) permet alors d'appliquer le th\'eor\`eme de Keller-Liverani \cite{keli} \cite{bal} 
\footnote{Les r\'esultats de \cite{keli} sont pr\'esent\'es 
avec une norme auxiliaire $|\cdot|$ sur $\cB$ v\'erifiant $|\cdot| \leq \|\cdot\|$, mais 
on peut montrer que ceux-ci subsistent lorsque $|\cdot|$ est remplac\'ee par une semi-norme, en l'occurrence ici $\nu(|\cdot|)$.}. 
Pour les assertions (a) et (b) nous appliquons ce th\'eor\`eme en $t_0=0$. Sous la condition (H1), celui-ci 
assure la d\'ecomposition (D) et l'assertion (a). 
Le point (b) est une cons\'equence de 
r\'esultats interm\'ediaires contenus dans \cite{keli} dont nous rappelons les principaux arguments. \\
Pr\'ecisons que le r\'eel $r$ de \cite{keli} [Th. 1] est choisi ici tel que 
$\frac{\ln \kappa - \ln r}{\ln\kappa} = \tau$, o\`u $\kappa$ est le r\'eel de (H4). 
Soit $\Gamma_1$ (resp. $\Gamma_0$) un cercle orient\'e de centre $z=1$, de rayon suffisamment petit (resp. de centre $z=0$, de 
rayon $\rho$ tel que $r < \rho <1$). 
D'apr\`es \cite{keli} [Th. 1], il existe $C>0$ telle que l'on ait pour $z\in\Gamma_0 \cup \Gamma_1$, 
$t\in J$, $f\in\cB$, 
$${\bf (**)} \ \ \ \ \  \ \ \ \ \ \ \nu\bigg(\bigg|(z-Q(t))^{-1}f -  (z-Q)^{-1}f\bigg|\bigg) \leq C\, |t|^\tau\, \|f\|.$$
Soit $\Pi_1(t) = \frac{1}{2i\pi}\ \int_{\Gamma_1}\ (z-Q(t))^{-1}\ dz$. Alors $\Pi_1(t)$ 
est un projecteur de rang 1 sur le sous-espace propre 
$\ker(Q(t)-\lambda(t))$ tel que $\Pi(0)f = \nu(f){\bf 1}$, et l'on peut d\'efinir 
$$v(t) = \langle \nu, \Pi(t){\bf 1}\rangle^{-1} \Pi(t){\bf 1}\ \ \ \ \mbox{et}\ \ \ \ \ \phi(t) = \Pi(t)^*\nu,$$
o\`u $\Pi(t)^*$ est l'op\'erateur adjoint de $\Pi(t)$ (on verra ci-dessous que $\langle \nu, \Pi(t){\bf 1}\rangle\neq 0$ 
pour $|t|$ petit 
de sorte que $v(t)$ est bien d\'efini). Le premier point de (b) est \'evident. \\ 
On a $\Pi(t){\bf 1} - {\bf 1} = \Pi(t){\bf 1} - \Pi(0){\bf 1} =  
\frac{1}{2i\pi} \int_{\Gamma_1}\ [(z-Q(t))^{-1}{\bf 1} -(z-Q)^{-1}{\bf 1}]\ dz$, d'o\`u 
$$\nu(|\Pi(t){\bf 1} - {\bf 1}|) \leq \frac{1}{2\pi} \int_{\Gamma_1}\ 
\nu\bigg(\bigg|(z-Q(t))^{-1}{\bf 1} - (z-Q)^{-1}{\bf 1}\bigg|\bigg)\ dz\ \leq\ C|\,|t|^\tau.$$
On en d\'eduit ais\'ement (b1), ainsi que (b2) gr\^ace \`a l'\'egalit\'e 
$\langle \phi(t),{\bf 1} \rangle = \langle \nu,\Pi(t){\bf 1} \rangle$. 
Enfin (b3) r\'esulte de l'\'egalit\'e 
$$N(t)^n{\bf 1} = N(t)^n{\bf 1} - N(0)^n{\bf 1} =  \frac{1}{2i\pi}\ \int_{\Gamma_0}\ z^n\, 
[(z-Q(t))^{-1}{\bf 1}  - (z-Q)^{-1}{\bf 1}] \ dz.$$
\noindent Il reste \`a prouver (c). Soit $t_0\in J$. Les deux premi\`eres majorations \'etablies au d\'ebut de la preuve  
et la propri\'et\'e (H4) montrent que $Q(t)$ v\'erifie au voisinage de $t_0$ les conditions du th\'eor\`eme de Keller-Liverani. 
Ce dernier, adjoint \`a la d\'ecomposition (D) \'ecrite en $t=t_0$, assure que, pour $|h|$ petit, $Q(t_0+h)$ 
admet une unique valeur propre dominante 
$z(t_0+h)$ qui tend vers $\lambda(t_0)$ quand $h\r0$. Mais par unicit\'e on a $z(t_0+h) = \lambda(t_0+h)$. \fdem 

\noindent En adaptant la preuve de \cite{ihp} [Lemme 4.2] \`a l'aide de la majoration du corolaire IV.1, 
nous allons maintenant pr\'eciser le comportement de $\lambda(u)$ quand $u\r0$. 

\noindent {\bf Lemme IV.2.} {\it Sous les hypoth\`eses ($\cH$) et $\sigma^2>0$, on a 
$\lambda(u) = 1-\frac{\sigma^2}{2}u^2 + O(|u|^{2+\tau})$ pour tout r\'eel $0<\tau<1$. \\ 
Si en outre $\langle \nu,|v(u)-{\bf 1}| \rangle = O(|u|)$, alors $\lambda(u) = 1-\frac{\sigma^2}{2}u^2 + O(u^3)$. }

\noindent{\it Preuve.} On suppose pour simplifier que $\sigma^2=1$. 
En utilisant le fait que $(X_n)_n$ est une cha\^{\i}ne de Markov (voir par exemple \cite{hulo} p. 23), puis 
le th\'eor\`eme IV, il est facile de voir que, pour tout $f\in\cB$, pour toute probabilit\'e initiale $\mu_0\in\cB'$, 
et pour $t\in J$, $n\geq1$, 
$$(***)\ \ \ \ \E\bigg[f(X_n)\, e^{itS_n}\bigg] = \langle \mu_0,Q(t)^nf \rangle = 
\lambda(t)^n \langle \phi(t),f \rangle\, \langle \mu_0,v(t)\rangle \ +\ 
\langle \mu_0, N(t)^n f \rangle.$$ 
Dans la suite on consid\`ere $t\in[-1,1]$ et un entier $n\geq N$, avec $N$ assez grand. Pour le moment $N$ est choisi tel que 
$\frac{t}{\sqrt n}\in J$ pour $n\geq N$. 
La formule ($***)$ appliqu\'ee en $\frac{t}{\sqrt n}$ avec $f=v(\frac{t}{ \sqrt n})$ et $\mu_0=\nu$ montre que  
$$\lambda(\frac{t}{ \sqrt n})^n =  \E[v(\frac{t}{ \sqrt n})(X_n)  e^{i  \frac{t}{ \sqrt n} S_n}].$$ 
En utilisant l'in\'egalit\'e triangulaire et l'invariance de $\nu$, il vient 
\begin{eqnarray*}
\bigg|\lambda(\mbox{$\frac{t}{ \sqrt n}$})^n - e^{-\frac{t^2}{2}}\bigg| &\leq&  
\E\bigg[\bigg|v(\mbox{$\frac{t}{ \sqrt n}$})(X_n) - 1 \bigg|\bigg]\ + 
\ \bigg|\E[e^{i\frac{t}{ \sqrt n}S_n}] - e^{-\frac{t^2}{2}}\bigg| \\
&=& \langle\nu,|v(\mbox{$\frac{t}{ \sqrt n}$}) - {\bf 1} |\rangle\ +\ 
|\E[e^{i\frac{t}{ \sqrt n}S_n}] - e^{-\frac{t^2}{2}}|.
\end{eqnarray*}
Du corollaire IV.1 et du point (b1) du th\'eor\`eme IV il vient pour $t\in[-1,1]$ et $n\geq N$, 
$$\bigg|\lambda(\mbox{$\frac{t}{ \sqrt n}$})^n - e^{-\frac{t^2}{2}}\bigg| \leq C\, |\frac{t}{\sqrt n}|^\tau.$$
Par ailleurs cette propri\'et\'e est satisfaite avec $\tau=1$ si $\langle \nu,|v(u)-{\bf 1}| \rangle = O(|u|)$. \\
En utilisant la fonction $\log$ complexe d\'efinie pour $z\in\C$ non nul par 
$\log z = \ln |z| + i\, \arg(z)$, avec $\arg(z)\in]-\pi,\pi]$, on d\'emontre qu'il existe 
une constante $C'>0$ telle que l'on ait pour $t\in[-1,1]$ et $n\geq N$, avec $N$ assez grand 
(voir les d\'etails dans \cite{ihp} p. 193, la continuit\'e de $\lambda(\cdot)$ sur $J$ est importante pour ce point), 
$$(L)\ \ \ \ \ \ \ \ \ \ \ \ \ \ \ \ \ \ \bigg|n\log \lambda(\mbox{$\frac{t}{\sqrt n}$}) + \frac{t^2}{2} \bigg| \leq 
C'\, |\frac{t}{\sqrt n}|^\tau.$$ 
En observant maintenant que $z-1 - \log z = (\log z)\, \alpha(z)$ avec $\alpha(z) = O(|z-1|)$, on a 
\begin{eqnarray*}
\bigg|n \bigg(\lambda(\mbox{$\frac{t}{ \sqrt n}$}) - 1\bigg) + \frac{t^2}{2} \bigg| &\leq& 
n \bigg|\lambda(\mbox{$\frac{t}{ \sqrt n}$}) - 1 - \log\lambda(\mbox{$\frac{t}{ \sqrt n}$})\bigg| \ +\   
\bigg|n \log\lambda(\mbox{$\frac{t}{ \sqrt n}$})  \ +\  \frac{t^2}{2} \bigg| \\ 
&\leq& \bigg|n\, \log\lambda(\mbox{$\frac{t}{ \sqrt n}$})\bigg|\ \bigg|\alpha(\lambda(\mbox{$\frac{t}{\sqrt n}$}))\bigg| 
+ C'\, |\frac{t}{\sqrt n}|^\tau 
\end{eqnarray*}
En utilisant \`a nouveau (L) et le fait que $|t|\leq 1$, 
on voit que $|n\log \lambda(\mbox{$\frac{t}{\sqrt n}$})|\leq \frac{1}{2} + C'$. \\
En outre le corollaire IV.1 montre que la suite $(\frac{S_n}{\sqrt n})_n$ converge en loi vers $\cN(0,\sigma^2)$. On d\'eduit alors  
de \cite{ihp} [Lemme 4.2] que $\lambda(u) = 1-\frac{u^2}{2} + o(u^2)$. 
Donc $\alpha(\lambda(u)) = O(\lambda(u) - 1) = O(u^2)$.  
Par cons\'equent il existe une constante $C''>0$ telle que l'on ait pour $t\in[-1,1]$ et $n\geq N$, avec $N$ assez grand, 
$$\bigg|n \bigg(\lambda(\mbox{$\frac{t}{ \sqrt n}$}) - 1\bigg) + \frac{t^2}{2} \bigg| \leq  
C''\, |\frac{t}{\sqrt n}|^\tau.$$
En divisant cette in\'egalit\'e par $t^2$ pour $\frac{1}{2} \leq |t| \leq 1$, on obtient 
$$\bigg|(\frac{t^2}{n})^{-1}\ \bigg(\lambda(\mbox{$\frac{t}{\sqrt n}$}) - 1\bigg) + 
\frac{1}{2}\bigg|\  \leq\  \frac{C''}{t^2}\, |\frac{t}{\sqrt n}|^\tau
\ \leq\  4\, C''\, |\frac{t}{\sqrt n}|^\tau.$$ 
Soit $u\in\R^*$, $|u|\leq \frac{1}{\sqrt{N}}$. Il existe clairement $n\geq N$ tel que 
$\frac{1}{2\sqrt n} \leq |u| \leq \frac{1}{\sqrt n}$. 
L'in\'egalit\'e pr\'ec\'edente appliqu\'ee avec $t = \sqrt n u$ 
montre que $|\frac{\lambda(u) - 1}{u^2} + \frac{1}{2}| \leq 4C''\, |u|^\tau$. Ceci d\'emontre le premier point du lemme. \\
Comme d\'ej\`a indiqu\'e, si $\langle \nu,|v(u)-{\bf 1}| \rangle = O(|u|)$, les arguments pr\'ec\'edents s'appliquent 
avec $\tau=1$. \fdem 

\noindent{\bf IV.3. D\'emonstration du th\'eor\`eme I }

\noindent Supposons pour simplifier que $\sigma^2=1$, et rappelons que, en vertu de l'in\'egalit\'e de Berry-Esseen,  
une vitesse en $n^{-\frac{\tau}{2}}$ ($\tau<1$) sera obtenue dans le t.l.c si l'on d\'emontre que, pour un certain $\alpha>0$, on a 
$$A_n = \int_{-\alpha\sqrt n}^{\alpha\sqrt n}
\bigg|\frac{\E[e^{it\frac{S_n}{\sqrt n}}] - e^{-\frac{t^2}{2}}}{t}\bigg|\ dt = 
O(n^{-\frac{\tau}{2}}).$$
Pour le moment on choisit $\alpha>0$ tel que $\frac{\alpha}{\sqrt n}\in J$, o\`u $J$ est l'intervalle du th\'eor\`eme IV. \\
En appliquant la formule ($***$) du paragraphe IV.2 avec $f={\bf 1}$, $\mu_0=\nu$, 
et en posant $L(u) = \langle \phi(u),{\bf 1} \rangle -1$, 
on a 
$\displaystyle \E[e^{it\frac{S_n}{\sqrt n}}] = \lambda(\frac{t}{\sqrt n})^n + \lambda(\frac{t}{\sqrt n})^n L(\frac{t}{\sqrt n})  \ +\ 
\langle \nu, N(\frac{t}{\sqrt n})^n {\bf 1} \rangle$. D'o\`u 
\begin{eqnarray*}
A_n & \leq &   
\int_{-\alpha\sqrt n}^{\alpha\sqrt n} \bigg|\frac{\lambda(\frac{t}{\sqrt n})^n - e^{-\frac{t^2}{2}}}{t}\bigg| \ dt+ 
\int_{-\alpha\sqrt n}^{\alpha\sqrt n} \bigg|\frac{\lambda(\frac{t}{\sqrt n})^n L(\frac{t}{\sqrt n})}{t}\bigg|\ dt + 
\int_{-\alpha\sqrt n}^{\alpha\sqrt n} \bigg|\frac{\langle \nu, N(\frac{t}{\sqrt n})^n {\bf 1} \rangle}{t}\bigg|\ dt \\
&=& I_n + J_n + K_n. 
\end{eqnarray*}
Du lemme IV.2, on d\'eduit que 
$|\lambda(u)| \leq  1-\frac{u^2}{4} \leq e^{-\frac{u^2}{4}}$ pour $|u|$ petit. 
Par cons\'equent on a pour $|\frac{t}{\sqrt n}|\leq \alpha$, avec $\alpha$ assez petit, 
$$|\lambda(\frac{t}{\sqrt n})| \leq e^{-\frac{t^2}{4n}},\ \ \mbox{d'o\`u}\ \ \  
|\lambda(\frac{t}{\sqrt n})|^n \leq e^{-\frac{t^2}{4}}.$$ 
Comme dans la preuve du th\'eor\`eme de Berry-Esseen, on \'ecrit 
$$\lambda(\frac{t}{\sqrt n})^n - e^{-\frac{t^2}{2}} = 
\bigg(\lambda(\frac{t}{\sqrt n}) - e^{-\frac{t^2}{2n}}\bigg) \, 
\sum_{k=0}^{n-1}  \lambda(\frac{t}{\sqrt n})^{n-k-1} e^{\frac{-kt^2}{2n}}.$$ 
Il existe $C,C' >$ tels que l'ait pour $\frac{|t|}{\sqrt n}$ assez petit, d'une part 
$|\lambda(\frac{t}{\sqrt n}) - e^{-\frac{t^2}{2n}}| \leq C |\frac{t}{\sqrt n}|^{2+\tau}$ (lemme IV.2), d'autre part 
$\displaystyle 
\sum_{k=0}^{n-1}  |\lambda(\frac{t}{\sqrt n})|^{n-k-1} e^{\frac{-kt^2}{2n}} \leq \sum_{k=0}^{n-1} e^{\frac{-t^2(n-k-1)}{4n}}\, 
^{\frac{-kt^2}{4n}} \leq C'n e^{-\frac{t^2}{4}}$, d'o\`u 
$$|\lambda(\frac{t}{\sqrt n})^n - e^{-\frac{t^2}{2}}| \leq CC'\, n^{-\frac{\tau}{2}}\, |t|^{2+\tau}\, 
e^{-\frac{t^2}{4}}.$$ 
La fonction $t\mapsto t^{1+\tau}\, e^{-\frac{t^2}{4}}$ \'etant int\'egrable sur $\R$, on a  $I_n = O(n^{-\frac{\tau}{2}})$. \\
En outre, de $|L(\frac{t}{\sqrt n})| \leq  C\, |\frac{t}{\sqrt n}|^{\tau}$ 
(Th. IV(b2)), et de l'int\'egrabilit\'e de $t\mapsto |t|^{\tau-1}\, e^{-\frac{t^2}{4}}$, 
on voit ais\'ement que $J_n = O(n^{-\frac{\tau}{2}})$. 
Enfin, du point (b3) du th\'eor\`eme IV, du fait que  $t\mapsto  t^{\tau-1}$ est int\'egrable sur $[-\alpha,\alpha]$ 
et que $\rho^n= O(n^{-\frac{\tau}{2}})$, on obtient $K_n = O(n^{-\frac{\tau}{2}})$. \fdem 

\noindent{\bf IV.4. D\'emonstration du th\'eor\`eme II}

\noindent On suppose ici que l'hypoth\`ese ($\widetilde{\cH}$) est satisfaite et que $\mu_0\in\cB'\cap\widetilde{\cB}'$, 
o\`u $\mu_0$ est la loi initiale. 
L'int\'egrale $A_n$ est d\'efinie comme en IV.3. 
En appliquant la formule ($***$) avec $f={\bf 1}$, on obtient en posant ici 
$L(u) = \langle \phi(u),{\bf 1} \rangle\, \langle \mu_0,v(u)\rangle - 1$, 
\begin{eqnarray*}
A_n & \leq &   
\int_{-\alpha\sqrt n}^{\alpha\sqrt n} \bigg|\frac{\lambda(\frac{t}{\sqrt n})^n - e^{-\frac{t^2}{2}}}{t}\bigg| \ dt+ 
\int_{-\alpha\sqrt n}^{\alpha\sqrt n} \bigg|\frac{\lambda(\frac{t}{\sqrt n})^n L(\frac{t}{\sqrt n}))}{t}\bigg| \ dt + 
\int_{-\alpha\sqrt n}^{\alpha\sqrt n} \bigg|\frac{\langle \mu_0, N(\frac{t}{\sqrt n})^n {\bf 1} \rangle}{t}\bigg|\ dt \\
&=& I_n + J_n + K_n, 
\end{eqnarray*}
Pour l'\'etude de $I_n$, $J_n$ et $K_n$, on consid\`ere 
des cercles orient\'es  $\Gamma_1$ et $\Gamma_0$ comme au $\S$ IV.2. En utilisant la derni\`ere condition de  ($\widetilde{\cH}$), 
il est facile de voir (cf. \cite{ihp} pp. 194) que, 
pour $f\in\cB$, $t\in J$ et $z\in\Gamma_0\cup\Gamma_1$, on a 
$$\|(z-Q(t))^{-1}f -  (z-Q)^{-1}f\|_{\sim} \leq C\, |t|\, \|f\|,$$
avec $C>0$ ind\'ependante de $f$, $t$ et $z$. Par int\'egration curviligne (voir les d\'efinitions de $v(t)$, $\phi(t)$, $N(t)^n$ au 
$\S$ IV.2), on obtient que \\
$\|v(u)-{\bf 1}\|_{\sim} \leq  C\, |u|$ \\
$|\langle \phi(u),{\bf 1} \rangle - 1| \leq  C\, |u|$ (utiliser le fait que $\nu\in \widetilde{\cB}'$). \\
$|\langle \mu_0, N(u)^n {\bf 1} \rangle| \leq  C\, \rho^n \, |u|$ (utiliser le fait que $\mu_0\in \widetilde{\cB}'$). \\
La premi\`ere in\'egalit\'e ci-dessus et le fait que $\widetilde{\cB}$ s'envoie contin\^ument dans $\L^1(\nu)$ montrent que 
$\langle \nu,|v(u)-{\bf 1}| \rangle = O(|u|)$. Du lemme IV.2, il vient que $\lambda(u) = 1-\frac{\sigma^2}{2}u^2 + O(u^3)$. 
Les arguments vus au $\S$ pr\'ec\'edent pour la majoration de $I_n$ s'appliquent alors avec $\tau=1$. 
Donc $I_n = O(n^{-\frac{1}{2}})$. En outre, des majorations ci-dessus, on peut d\'eduire que $L(u) = O(|u|)$, et l'on a vu que 
$|\lambda(\frac{t}{\sqrt n})|^n \leq e^{-\frac{t^2}{4}}$ pour $|\frac{t}{\sqrt n}|$ assez petit ; on en d\'eduit que 
$J_n = O(n^{-\frac{1}{2}})$. Enfin on a clairement $K_n = O(\rho^n) = O(n^{-\frac{1}{2}})$. \fdem

\end{document}